\documentclass[11pt]{article}%
\usepackage{amsfonts}
\usepackage{amsmath}
\usepackage{amssymb}
\usepackage{graphicx}%
\providecommand{\U}[1]{\protect\rule{.1in}{.1in}}
\everymath{\displaystyle}
\newtheorem{theorem}{Theorem}

\newtheorem{definition}[theorem]{Definition}
\newtheorem{example}[theorem]{Example}

\newtheorem{proposition}[theorem]{Proposition}

\newenvironment{proof}[1][Proof]{\noindent\textbf{#1.} }{\ \rule{0.5em}{0.5em}}
\begin{document}

\date{}
\title{Orthogonality of polar Legendre polynomials and approximation}
\author{Abdelhamid Rehouma\\Department of mathematics Faculty of exact sciences\\University Hama Lakhdar, Eloued Algeria.\\URL: https://sites.google.com/view/mathsrehoumablog/accueil\\E-mail : rehoumaths@gmail.com}
\maketitle

\begin{abstract}
Let $\left\{  Q_{n}\left(  x\right)  \right\}  $ be a system of integral
Legendre polynomials of degree exactly $n$,and let $\left\{  P_{n}\left(
x\right)  \right\}  $ be polar polynomials primitives of integral Legendre
polynomials .We derive some identities and relations and extremal problems and
minimization involving of polar integral Legendre polynomials.

\end{abstract}

\bigskip\textbf{Keywords:}{\scriptsize \ Legendre \ polynomials, PIPCIR
polynomials,\ Three-term rrecursive relation, Summations;Polar Legendre
polynomials. \vspace{4mm} \newline}

\noindent\textbf{2000. MSC:}{\footnotesize \ 42C05, 33C45.}

\section{Mathematical basis}

We restricted our attention to a polynomial with the first and last roots at
$x=\pm1$ ,given by
\begin{equation}
Q_{n}\left(  x\right)  =\left(  x^{2}-1\right)  q_{n-2}\left(  x\right)
\text{ \ \ \ \ ,}n\geq2 \label{Defi1}%
\end{equation}
Let us call a polynomial whose inflection points coincide with their interior
roots in a shorter way : Pipcir .It will be shown that the zeros of these
polynomials are all real, distinct, and they lie in the interval $\left[
-1\text{ \ \ \ \ }1\right]  .$The requirement all inflection points to
coincide with all roots of $Q_{n}\left(  x\right)  $ except $\pm1$ yields:%
\[
Q_{n}^{\prime\prime}\left(  x\right)  =-n\left(  n-1\right)  q_{n-2}\left(
x\right)
\]
or%
\begin{equation}
\left(  1-x^{2}\right)  Q_{n}^{\prime\prime}\left(  x\right)  +n\left(
n-1\right)  Q_{n}\left(  x\right)  =0 \label{Diff2}%
\end{equation}
Let us differentiate the equation \eqref{Diff2}%
\begin{equation}
\left(  1-x^{2}\right)  Q_{n}^{\prime\prime\prime}\left(  x\right)
-2xQ_{n}^{\prime\prime}\left(  x\right)  +n\left(  n-1\right)  Q_{n}^{\prime
}\left(  x\right)  =0 \label{Diff3}%
\end{equation}

We have now well-known Legendre's differential equation whose bounded on
$\left[  -1\text{ \ \ \ \ }1\right]  $ solutions are known as Legendre
polynomials:$y_{n}=L_{n-1}\left(  x\right)  $ ,$n\geq1.$One can find
properties of these polynomials in \cite{Abram},\cite{Belinsk},\cite{Pijeir}%
,they are normalized so that $L_{n}\left(  1\right)  =1$ for all $n$ .If%
\begin{equation}
Q_{n}\left(  x\right)  =-%
{\displaystyle\int\limits_{x}^{1}}
L_{n-1}\left(  t\right)  dt\text{ \ \ \ \ \ }-1\leq x\leq1 \label{ExpQn}%
\end{equation}
then $Q_{n}^{\prime}\left(  x\right)  =L_{n-1}\left(  x\right)  $ and
$Q_{n}^{\prime\prime}\left(  x\right)  =L_{n-1}^{\prime}\left(  x\right)  $
.We see that polynomials $Q_{n}\left(  x\right)  $ defined by \eqref{ExpQn}
satisfy the equation \eqref{Diff2}.Thus,%
\begin{equation}
Q_{n}\left(  1\right)  =Q_{n}\left(  -1\right)  =0\text{ \ \ \ \ ,}n\geq2
\label{Qqn1}%
\end{equation}

The Legendre Polynomial,$L_{n}\left(  x\right)  $ saisfies,\cite{Abram}%
,\cite{Belinsk},\cite{Pijeir} :%
\begin{equation}
L_{n}\left(  x\right)  =\frac{1}{2^{n}n!}\left(  \left(  x^{2}-1\right)
^{n}\right)  ^{\left(  n\right)  } \label{DifLn}%
\end{equation}

The Legendre polynomials, denoted by $L_{k}\left(  x\right)  $, are the
orthogonal polynomials with $\omega\left(  x\right)  =1.$The three-term
recurrence relation for the Legendre polynomials reads,\cite{Abram}%
,\cite{Belinsk},\cite{Pijeir}%
\[
L_{0}\left(  x\right)  =1,L_{1}\left(  x\right)  =x,
\]
and%
\[
\left(  n+1\right)  L_{n+1}\left(  x\right)  =\left(  2n+1\right)
xL_{n}\left(  x\right)  -nL_{n-1}\left(  x\right)  \text{ \ \ \ \ \ }%
n=1,2,3..
\]
They are normalized so that $L_{n}\left(  1\right)  =1$ for all $n.$
yields,\cite{Abram},\cite{Belinsk},\cite{Pijeir}:%
\[
\left(  \left(  1-x^{2}\right)  L_{n}^{\prime}\left(  x\right)  \right)
^{\prime}+n\left(  n+1\right)  L_{n}\left(  x\right)  =0
\]
and%
\begin{equation}%
{\displaystyle\int\limits_{-1}^{1}}
L_{n}\left(  x\right)  L_{m}\left(  x\right)  dx=\frac{2}{2n+1}\delta
_{n,m},\text{ \ \ \ \ }n,m=1,2,3.. \label{orthLn}%
\end{equation}
and%
\[%
{\displaystyle\int\limits_{-1}^{1}}
L_{n}^{2}\left(  x\right)  =\frac{2n-1}{2n+1}%
{\displaystyle\int\limits_{-1}^{1}}
L_{n-1}^{2}\left(  x\right)  \text{ \ \ \ \ }n=1,2,3..
\]
We also derive that,\cite{Abram},\cite{Belinsk},\cite{Pijeir}%
\[%
{\displaystyle\int\limits_{-1}^{x}}
L_{n}\left(  t\right)  dt=\frac{1}{2n+1}\left(  L_{n+1}\left(  x\right)
-L_{n-1}\left(  x\right)  \right)
\]
We derive from above a recursive relation for computing the derivatives of the
Legendre polynomials,\cite{Abram},\cite{Belinsk},\cite{Pijeir}:
\[
L_{n}\left(  x\right)  =\frac{1}{2n+1}\left(  L_{n+1}^{\prime}\left(
x\right)  -L_{n-1}^{\prime}\left(  x\right)  \right)
\]
We also derive that,\cite{Abram}%
\begin{equation}
L_{n}\left(  \pm1\right)  =\left(  \pm1\right)  ^{n} \label{Lnat1}%
\end{equation}%
\[
L_{n}^{\prime}\left(  \pm1\right)  =\frac{1}{2}\left(  \pm1\right)
^{n-1}n\left(  n+1\right)
\]%
\[
L_{n}^{\prime\prime}\left(  \pm1\right)  =\left(  \pm1\right)  ^{n}\left(
n-1\right)  n\left(  n+1\right)  \left(  n+2\right)  /8
\]
with%
\begin{equation}
Q_{n}^{\prime}\left(  1\right)  =1 \label{Qnderiv1}%
\end{equation}%
\begin{equation}
Q_{n}^{\prime\prime}\left(  1\right)  )=\frac{1}{2}\left(  \pm1\right)
^{n-1}n\left(  n-1\right)  \label{Secondderiqat1}%
\end{equation}
We also derive that%
\begin{equation}
L_{n}\left(  x\right)  =\frac{1}{2^{n}}\sum_{k=0}^{n}\left(  C_{k}^{n}\right)
^{2}\left(  x-1\right)  ^{n-k}\left(  x+1\right)  ^{k} \label{expp}%
\end{equation}
becomes%
\begin{equation}
L_{n}\left(  0\right)  =\frac{1}{2^{n}}\sum_{k=0}^{n}\left(  -1\right)
^{n-k}\left(  C_{k}^{n}\right)  ^{2} \label{Lnat0}%
\end{equation}
and%
\begin{equation}
L_{n}^{\prime}\left(  0\right)  =\frac{1}{2^{n}}\sum_{k=0}^{n}\left(
-1\right)  ^{n-k}\left(  -n+2k\right)  \left(  C_{k}^{n}\right)  ^{2}
\label{Lnderivat0}%
\end{equation}
Explicit formula for $Q_{n}\left(  x\right)  $is the following, \cite{Belinsk}%
:%
\begin{equation}
Q_{n}\left(  x\right)  =%
{\displaystyle\sum\limits_{k=0}^{n}}
\frac{\left(  -1\right)  ^{k}\left(  2n-2k-3\right)  !!}{\left(  2k\right)
!!\left(  n-2k\right)  !!}x^{n-2k} \label{expliv}%
\end{equation}
and%
\begin{equation}
Q_{n}\left(  0\right)  =\frac{\left(  -1\right)  ^{\frac{n-2}{2}}\left(
n-3\right)  !!}{n!!} \label{Qnat0}%
\end{equation}
The Rodrigues formula for the $\left\{  Q_{n}\right\}  _{n=2,3,4.....}$
orthogonal polynomials is well known as the following ,\cite{Belinsk}%
\begin{equation}
Q_{n}\left(  x\right)  =\frac{x^{2}-1}{2^{n-1}n!\left(  n-1\right)  }\left[
\left(  x^{2}-1\right)  ^{n-1}\right]  ^{\left(  n\right)  } \label{Rodrigues}%
\end{equation}
Now we have two expressions for $Q_{n}\left(  x\right)  $; equating them, we
obtain the formula%
\begin{equation}
\left(  x^{2}-1\right)  \left[  \left(  x^{2}-1\right)  ^{n-1}\right]
^{\left(  n\right)  }=n\left(  n-1\right)  \left[  \left(  x^{2}-1\right)
^{n-1}\right]  ^{\left(  n-2\right)  } \label{Second}%
\end{equation}
As is well known \cite{Abram},we note that%
\[%
{\displaystyle\int\limits_{-1}^{1}}
uv^{\left(  n\right)  }=\left\{
{\displaystyle\sum\limits_{k=1}^{n}}
\left(  -1\right)  ^{k-1}u^{\left(  k-1\right)  }v^{\left(  n-k\right)
}\right\}  _{-1}^{1}+\left(  -1\right)  ^{n}%
{\displaystyle\int\limits_{-1}^{1}}
u^{\left(  n\right)  }v
\]
We derive from above a recursive relation for computing the derivatives of the
\textbf{Legendr}e polynomials,\cite{Abram}:
\[
L_{n}\left(  x\right)  =\frac{1}{2n+1}\left(  L_{n+1}^{\prime}\left(
x\right)  -L_{n-1}^{\prime}\left(  x\right)  \right)
\]
reduces to%
\begin{equation}
Q_{n}^{\prime}\left(  x\right)  =\frac{1}{2n-1}\left(  Q_{n+1}^{\prime\prime
}\left(  x\right)  -Q_{n-1}^{\prime\prime}\left(  x\right)  \right)
\label{Pipcirs2}%
\end{equation}
Integrating both sides \ of \eqref{Pipcirs2} yields%
\begin{equation}
\int Q_{n}\left(  x\right)  dx=\frac{1}{2n-1}\left(  Q_{n+1}\left(  x\right)
-Q_{n-1}\left(  x\right)  \right)  \label{Pipcirs3}%
\end{equation}

\subsection{ Polar Legendre polynomials and orthogonality}

\begin{definition}
Let $\ P_{n}$ define as apolynomial of degree $n$ such that%
\begin{equation}
-\left(  n+1\right)  \int_{x}^{1}L_{n}\left(  z\right)  dz=\left(  x-1\right)
P_{n}\left(  x\right)  \label{polar}%
\end{equation}

normalised by
\begin{equation}
\left[  \left(  x-1\right)  P_{n}\left(  x\right)  \right]  _{x=1}=0
\label{condition}%
\end{equation}

such that%
\begin{equation}
\left(  n+1\right)  L_{n}\left(  x\right)  =\left[  \left(  x-1\right)
P_{n}\left(  x\right)  \right]  ^{\prime}=\left(  x-1\right)  P_{n}^{\prime
}\left(  x\right)  +P_{n}\left(  x\right)  \label{condpolar}%
\end{equation}

$P_{n}$ is called the $n$-th polar\textbf{\ Legendre }polynomial.
Obviously,$P_{n}$ is a polynomial of degree $n$, This type of polar
\textbf{Legendre} polynomials was introduced and studied initially in

$[1]$.Obviously the following calculus shows that the pole of \ $P_{n}\left(
x\right)  $ do not have to be irregular.%
\begin{equation}
\lim_{z\longrightarrow1}P_{n}\left(  x\right)  =\left(  n+1\right)
\lim_{x\longrightarrow1}\frac{-\int_{x}^{1}L_{n}\left(  z\right)  dz}%
{x-1}=\left(  n+1\right)  L_{n}\left(  1\right)
\text{\ \ \ \ \ \ \ \ \ \ \ \ } \label{pole}%
\end{equation}

it is appear that%
\begin{equation}
\left(  n+1\right)  Q_{n+1}\left(  x\right)  =\left(  x-1\right)  P_{n}\left(
x\right)  \label{QRP}%
\end{equation}

From \eqref{Rodrigues} it follows that
\begin{equation}
P_{n}\left(  x\right)  =\frac{x+1}{2^{n}n!n}\left[  \left(  x^{2}-1\right)
^{n}\right]  ^{\left(  n+1\right)  } \label{PRodrigues}%
\end{equation}

\end{definition}

\section{Identities and relations involving polar Legendre polynomials}

\begin{proposition}
The polynomials $P_{n}$ satisfy the following linear differential equation :%
\begin{equation}
\left(  x^{2}-1\right)  P_{n}^{\prime\prime}\left(  x\right)  +2\left(
x+1\right)  P_{n}^{\prime}\left(  x\right)  -n\left(  n+1\right)  P_{n}\left(
x\right)  =0 \label{Polardiffequat}%
\end{equation}

and%
\begin{equation}
P_{n}\left(  0\right)  =\frac{\left(  -1\right)  ^{\frac{n}{2}}\left(
n+1\right)  \left(  n-3\right)  !!}{n!!} \label{Pnat0}%
\end{equation}

and%
\begin{equation}
P_{n}\left(  1\right)  =n+1 \label{Pnat1}%
\end{equation}

with%
\begin{equation}
P_{n}^{\prime}\left(  1\right)  =\frac{n\left(  n^{2}-1\right)  }{4}
\label{derivpnat1}%
\end{equation}

\end{proposition}

\begin{proof}
Let us differentiate\eqref{QRP} the equation with respect to $x$:%
\begin{equation}
\left(  n+1\right)  Q_{n+1}^{\prime}\left(  x\right)  =P_{n}\left(  x\right)
+\left(  x-1\right)  P_{n}^{\prime}\left(  x\right)  \label{derivPQn}%
\end{equation}

from which it follows%
\begin{equation}
\left(  n+1\right)  Q_{n+1}^{\prime\prime}\left(  x\right)  =2P_{n}^{\prime
}\left(  x\right)  +\left(  x-1\right)  P_{n}^{\prime\prime}\left(  x\right)
\label{SeconDerivPQn}%
\end{equation}

applying \eqref{Diff2},\eqref{Diff3},\eqref{QRP}, we can deduce after several
computations that%
\[
\left(  1-x^{2}\right)  \left(  2P_{n}^{\prime}\left(  x\right)  +\left(
x-1\right)  P_{n}^{\prime\prime}\left(  x\right)  \right)  +n\left(
n+1\right)  \left(  x-1\right)  P_{n}\left(  x\right)  =0
\]

we have%
\[
\left(  x^{2}-1\right)  P_{n}^{\prime\prime}\left(  x\right)  +2\left(
x+1\right)  P_{n}^{\prime}\left(  x\right)  -n\left(  n+1\right)  P_{n}\left(
x\right)  =0
\]

and \eqref{Polardiffequat} is proved.We can put%
\[
\left(  n+1\right)  Q_{n}\left(  0\right)  =-P_{n}\left(  0\right)
\]

using \eqref{Qnat0},we get:%
\[
P_{n}\left(  0\right)  =-\frac{\left(  -1\right)  ^{\frac{n-2}{2}}\left(
n+1\right)  \left(  n-3\right)  !!}{n!!}%
\]

from \eqref{QRP} \eqref{Qnderiv1}we deduce the following result directly:%
\begin{equation}
P_{n}\left(  1\right)  =\left(  n+1\right)  Lim_{x\longrightarrow1}\frac
{Q_{n}\left(  x\right)  }{x-1}=\left(  n+1\right)  Q_{n}^{\prime}\left(
1\right)  =n+1 \label{Pnt1}%
\end{equation}

Using \eqref{Polardiffequat}%
\[
4P_{n}^{\prime}\left(  1\right)  -n\left(  n-1\right)  P_{n}\left(  1\right)
=0
\]

By \eqref{Pnat1} we deduce
\[
P_{n}^{\prime}\left(  1\right)  =\frac{n\left(  n^{2}-1\right)  }{4}%
\]

and the proposition is proved.
\end{proof}

\section{Main results}

\subsection{Ortogonality of polar Legendre polynomials}

You may see examples of polynomials $Q_{n}\left(  x\right)  ,$see
\cite{Belinsk}
\[
Q_{2}\left(  x\right)  =\frac{1}{2}\left(  x^{2}-1\right)
\]%
\[
Q_{3}\left(  x\right)  =\frac{1}{2}\left(  x^{3}-x\right)
\]%
\[
Q_{4}\left(  x\right)  =\frac{1}{8}\left(  5x^{4}-6x^{2}+1\right)
\]%
\[
Q_{5}\left(  x\right)  =\frac{1}{8}\left(  7x^{5}-10x^{3}+3x\right)
\]%
\[
Q_{6}\left(  x\right)  =\frac{1}{16}\left(  21x^{6}-35x^{4}+15x^{2}-1\right)
,
\]
First we prove that the functions $Q_{n}\left(  x\right)  $ and $Q_{m}\left(
x\right)  $) ($n\neq m$) are orthogonal over $\left[  -1\text{ \ \ \ }%
1\right]  ,$with respect to the weight function $w\left(  x\right)  =\dfrac
{1}{1-x^{2}}.$

\begin{theorem}
We have
\end{theorem}

\begin{equation}%
{\displaystyle\int\limits_{-1}^{1}}
\dfrac{Q_{n}\left(  x\right)  Q_{m}\left(  x\right)  }{1-x^{2}}%
dx=0,\ \ \ \ \ (n\neq m),\ \ \ n,m=1,2,3 \label{Orthogo}%
\end{equation}

and%
\begin{equation}
\left\Vert Q_{n}\right\Vert ^{2}=%
{\displaystyle\int\limits_{-1}^{1}}
\dfrac{Q_{n}^{2}\left(  x\right)  }{1-x^{2}}dx=\frac{2}{n\left(  n-1\right)
\left(  2n-1\right)  },\ \ \ \ \ \ \ \ \ \ \ n=2,3,4..... \label{NormQn}%
\end{equation}

but because $Q_{n}\left(  x\right)  =0,Q_{n}\left(  -1\right)  =0$, all
integrals \eqref{Orthogo},are proper.

\begin{proof}
Now formulas \eqref{Diff2}and\eqref{Rodrigues}it follows that, for
$k=0,1,2,3.....n$%
\[
\dfrac{Q_{n}\left(  x\right)  x^{k}}{1-x^{2}}=\frac{-1}{n\left(  n-1\right)
}Q_{n}^{\prime\prime}\left(  x\right)  x^{k}=\frac{-1}{2^{n-1}n!n\left(
n-1\right)  ^{2}}\left(  \left(  x^{2}-1\right)  ^{n-1}\right)  ^{\left(
n+2\right)  }x^{k}%
\]

we obtain relation%
\[%
{\displaystyle\int\limits_{-1}^{1}}
\dfrac{Q_{n}\left(  x\right)  x^{k}}{1-x^{2}}dx=\frac{-1}{2^{n-1}n!n\left(
n-1\right)  ^{2}}%
{\displaystyle\int\limits_{-1}^{1}}
\left(  \left(  x^{2}-1\right)  ^{n-1}\right)  ^{\left(  n+2\right)  }x^{k}dx
\]%
\[
=\frac{-1}{2^{n-1}n!n\left(  n-1\right)  ^{2}}\left[  \left(  \left(
x^{2}-1\right)  ^{n-1}\right)  ^{\left(  n+1\right)  }x^{k}\right]
_{x=-1}^{x=1}+\frac{k}{2^{n-1}n!n\left(  n-1\right)  ^{2}}%
{\displaystyle\int\limits_{-1}^{1}}
\left(  \left(  x^{2}-1\right)  ^{n-1}\right)  ^{\left(  n+1\right)  }%
x^{k-1}dx
\]%
\[
=-\frac{k\left(  k-1\right)  }{2^{n-1}n!n\left(  n-1\right)  ^{2}}%
{\displaystyle\int\limits_{-1}^{1}}
\left(  \left(  x^{2}-1\right)  ^{n-1}\right)  ^{\left(  n\right)  }x^{k-2}dx
\]%
\[
................................................................................................
\]%
\[
=\pm\frac{k!}{2^{n-1}n!n\left(  n-1\right)  ^{2}}%
{\displaystyle\int\limits_{-1}^{1}}
\left(  \left(  x^{2}-1\right)  ^{n-1}\right)  ^{\left(  n-k+2\right)  }dx
\]%
\[
=\pm\frac{k!}{2^{n-1}n!n\left(  n-1\right)  ^{2}}\left[  \left(  \left(
x^{2}-1\right)  ^{n-1}\right)  ^{\left(  n-k+1\right)  }\right]  _{x=-1}%
^{x=1}=0
\]

Thus property \eqref{Orthogo} is proved .To prove \eqref{NormQn},we can see%
\[%
{\displaystyle\int\limits_{-1}^{1}}
\dfrac{Q_{n}\left(  x\right)  x^{n}}{1-x^{2}}dx=\pm\frac{1}{2^{n-1}\left(
n-1\right)  }%
{\displaystyle\int\limits_{-1}^{1}}
\left(  x^{2}-1\right)  ^{2\left(  n-1\right)  }dx
\]

In fact%
\[%
{\displaystyle\int\limits_{0}^{\frac{\pi}{2}}}
\sin^{2n-2}xdx=\frac{\left(  2n-2\right)  !}{4^{n-1}\left(  \left(
n-1\right)  !\right)  ^{2}}\frac{\pi}{2}%
\]

for mor details see,\cite{Belinsk} and the theorem is proved.
\end{proof}

Second we prove that the functions\ polar Legrndre polynomials,$P_{n}\left(
x\right)  $ and $P_{m}\left(  x\right)  $ ($n\neq m$) are orthogonal over
$\left[  -1\text{ \ \ \ }1\right]  ,$with respect to the weight function
$w\left(  x\right)  =\dfrac{1-x}{1+x}.$

\begin{theorem}
We have
\end{theorem}

\begin{equation}%
{\displaystyle\int\limits_{-1}^{1}}
P_{n}\left(  x\right)  P_{m}\left(  x\right)  \dfrac{1-x}{1+x}%
dx=0\ \ \ \ \ \ ,m\neq n,m,n=0,1,2. \label{Orthog}%
\end{equation}

and%

\begin{equation}%
{\displaystyle\int\limits_{-1}^{1}}
P_{n}^{2}\left(  x\right)  \dfrac{1-x}{1+x}dx=\frac{2\left(  n+1\right)  ^{2}%
}{n\left(  n-1\right)  \left(  2n-1\right)  }\text{ \ \ \ \ \ ,}n=2,3,4.......
\label{NormPn}%
\end{equation}

\begin{proof}
Combining the formulas \eqref{QRP}, \eqref{Orthog},\eqref{NormQn}, for $n\neq
m,n,m=2,3,4.....$%
\[%
{\displaystyle\int\limits_{-1}^{1}}
P_{n}\left(  x\right)  P_{m}\left(  x\right)  \dfrac{1-x}{1+x}dx=%
{\displaystyle\int\limits_{-1}^{1}}
P_{n}\left(  x\right)  P_{m}\left(  x\right)  \dfrac{\left(  x-1\right)  ^{2}%
}{1-x^{2}}dx
\]%
\[
=\left(  n+1\right)  \left(  m+1\right)
{\displaystyle\int\limits_{-1}^{1}}
Q_{n}\left(  x\right)  Q_{m}\left(  x\right)  \dfrac{dx}{1-x^{2}}=0
\]

and%
\[%
{\displaystyle\int\limits_{-1}^{1}}
P_{n}^{2}\left(  x\right)  \dfrac{1-x}{1+x}dx=%
{\displaystyle\int\limits_{-1}^{1}}
P_{n}^{2}\left(  x\right)  \dfrac{\left(  x-1\right)  ^{2}}{1-x^{2}}dx=\left(
n+1\right)  ^{2}%
{\displaystyle\int\limits_{-1}^{1}}
Q_{n}^{2}\left(  x\right)  \dfrac{dx}{1-x^{2}}%
\]

i-e%
\[
\left\Vert P_{n}\right\Vert ^{2}=\frac{2\left(  n+1\right)  ^{2}}{n\left(
n-1\right)  \left(  2n-1\right)  }\ \
\]

and the theorem is proved.
\end{proof}

\subsection{Kernels polynomials and extremal problem and minimization}

The $n$-th $Q$-kernel is given by, \cite{Reh}, \cite{Szeg}
\begin{equation}
K_{n}\left(  x,y\right)  =%
{\displaystyle\sum\limits_{k=0}^{n}}
\frac{P_{k}(x)P_{k}(y)}{\left\Vert P_{k}\right\Vert ^{2}}. \label{Kern}%
\end{equation}
satisfies the \textbf{Christoffel-Darboux} formula, \cite{Szeg},
\cite{Abramo},\cite{Reh}%
\begin{equation}
K_{n}\left(  x,y\right)  =\frac{1}{\left\Vert P_{n}\right\Vert ^{2}}%
\frac{P_{n+1}(x)P_{n}(y)-P_{n+1}(y)P_{n}(x)}{x-y},\ \ \ \ \ \ x\neq y
\label{CDS11}%
\end{equation}
and for $x=y$ one has
\begin{equation}
K_{n}\left(  x,x\right)  =\frac{1}{\left\Vert P_{n}\right\Vert ^{2}}\left(
P_{n+1}^{\prime}(x)P_{n}(x)-P_{n+1}(x)P_{n}^{\prime}(x)\right)  . \label{ABC}%
\end{equation}
$K_{n}$ has the reproducing kernel property \cite{Szeg}, \cite{Abramo}%
,\cite{Reh}:
\begin{equation}
f\left(  x\right)  =\int\limits_{-1}^{1}K_{n}\left(  x,t\right)  f\left(
t\right)  \dfrac{1-t}{1+t}dt \label{Reprkernel}%
\end{equation}
According to \eqref{Kern},%
\[
K_{n}\left(  x,0\right)  =%
{\displaystyle\sum\limits_{k=0}^{n}}
\frac{P_{k}(x)P_{k}(0)}{\left\Vert P_{k}\right\Vert ^{2}}%
\]
Combining the formulas \eqref{Pnat0}, \eqref{NormPn}%
\[
K_{n}\left(  x,0\right)  =%
{\displaystyle\sum\limits_{k=0}^{n}}
\left(  -1\right)  ^{\frac{k}{2}}\frac{k\left(  k-1\right)  \left(
2k-1\right)  \left(  k-3\right)  !!}{2\left(  k+1\right)  k!!}P_{k}(x)
\]
Hence%
\begin{equation}
K_{n}\left(  0,0\right)  =%
{\displaystyle\sum\limits_{k=0}^{n}}
\left(  -1\right)  ^{k}\frac{k\left(  k-1\right)  \left(  2k-1\right)  \left(
\left(  k-3\right)  !!\right)  ^{2}}{2\left(  k!!\right)  ^{2}} \label{Knat00}%
\end{equation}
The sequence $\left(  K_{n}\left(  x,0\right)  \right)  _{n=0}^{\infty}$ is
orthogonal with respect to the weight function
\[
t\left(  x\right)  =\dfrac{x\left(  1-x\right)  }{1+x}%
\]
for \ $-1\leq x\leq1$, i. e.%
\[
\int\limits_{-1}^{1}K_{n}\left(  x,0\right)  K_{m}\left(  x,0\right)
\dfrac{x\left(  1-x\right)  }{1+x}dx=0,\ n\neq m.
\]
According to \eqref{ABC}
\[
K_{n}\left(  0,0\right)  =\frac{1}{\left\Vert P_{n}\right\Vert ^{2}}\left(
P_{n+1}^{\prime}(0)P_{n}(0)-P_{n+1}(0)P_{n}^{\prime}(0)\right)
\]
To compute $P_{n}^{\prime}\left(  0\right)  $ using
\eqref{condpolar},\eqref{Pnat0},\eqref{Lnat0},\eqref{condpolar},%
\begin{equation}
P_{n}^{\prime}\left(  0\right)  =-\left(  n+1\right)  L_{n}\left(  0\right)
+P_{n}\left(  0\right)  \label{DerivPnat0}%
\end{equation}
where%
\[
P_{n}\left(  0\right)  =\frac{\left(  -1\right)  ^{\frac{n}{2}}\left(
n+1\right)  \left(  n-3\right)  !!}{n!!}%
\]
and%
\[
L_{n}\left(  0\right)  =\frac{1}{2^{n}}\sum_{k=0}^{n}\left(  -1\right)
^{n-k}\left(  C_{k}^{n}\right)  ^{2}%
\]
it follows that,%
\begin{equation}
P_{n}^{\prime}\left(  0\right)  =\frac{\left(  -1\right)  ^{\frac{n}{2}%
}\left(  n+1\right)  \left(  n-3\right)  !!}{n!!}-\frac{\left(  n+1\right)
}{2^{n}}\sum_{k=0}^{n}\left(  -1\right)  ^{n-k}\left(  C_{k}^{n}\right)  ^{2}
\label{DerivPnat00}%
\end{equation}
Using \eqref{Pnat0},\eqref{DerivPnat0},\eqref{NormPn},\eqref{DerivPnat00}we
deduce that%
\[
K_{n}\left(  0,0\right)  =
\]%
\[
=%
{\displaystyle\sum\limits_{k=0}^{n}}
\left(  -1\right)  ^{k}\frac{k\left(  k-1\right)  \left(  2k-1\right)  \left(
\left(  k-3\right)  !!\right)  ^{2}}{2\left(  k!!\right)  ^{2}}%
\]%
\[
\left(  -1\right)  ^{\frac{2n+1}{2}}\frac{n\left(  n-1\right)  \left(
2n-1\right)  \left(  n+1\right)  \left(  n+2\right)  \left(  n-2\right)
!!\left(  n-3\right)  !!}{2\left(  n+1\right)  ^{2}\left(  n+1\right)  !!n!!}%
\]%
\[
+\left(  -1\right)  ^{\frac{2n+3}{2}}\frac{n\left(  n-1\right)  \left(
2n-1\right)  \left(  n+2\right)  \left(  n+1\right)  \left(  n-3\right)
!!\left(  n-2\right)  !!}{2\left(  n+1\right)  ^{2}n!!\left(  n+1\right)  !!}%
\]%
\[
+\left(  -1\right)  ^{\frac{n+1}{2}}\frac{n\left(  n-1\right)  \left(
2n-1\right)  \left(  n+1\right)  \left(  n+2\right)  \left(  n-3\right)
!!}{2^{n+2}\left(  n+1\right)  ^{2}n!!}\sum_{k=0}^{n+1}\left(  -1\right)
^{n-k+1}\left(  C_{k}^{n+1}\right)  ^{2}%
\]%
\[
+\left(  -1\right)  ^{\frac{n+1}{2}}\frac{n\left(  n-1\right)  \left(
2n-1\right)  \left(  n+2\right)  \left(  n-2\right)  !!}{2^{n+1}\left(
n+1\right)  \left(  n+1\right)  !!}\sum_{k=0}^{n}\left(  -1\right)
^{n-k}\left(  C_{k}^{n}\right)  ^{2}%
\]

\subsection{Extremal problem and minimization}

Let $x\rightarrow w(x)=$ $\dfrac{1-x}{1+x}$ be a nonnegative function on the
interval $[-1,1]$ such that
\[
\int\limits_{-1}^{1}x^{r}w\left(  x\right)  dx
\]
exists for $r\geq0,$and consider the definite integral of the form%
\begin{equation}
I_{n}=\int\limits_{-1}^{1}f_{n}^{2}\left(  x\right)  \dfrac{1-x}{1+x}dx
\label{ExtrPbm}%
\end{equation}
where $f_{n}\left(  x\right)  $ is any real polynomial of degree n such that
$f_{n}\left(  1\right)  =1$.The problem to be solved is to determine the
polynomial $x\longrightarrow f_{n}\left(  x\right)  $ of order $n$ \ which
minimizes the integral \eqref{ExtrPbm} Since the integrand is non negative for
any value of $x\in\lbrack-1,1]$ such a minimum value does exist.

Using standard minimization technique \cite{Reh}, \cite{Szeg} ,and starting
from%
\[
\varphi\left(  a_{0},a_{1},....a_{n},\beta\right)  =\int\limits_{-1}%
^{1}\left(  \sum_{k=0}^{n}a_{k}P_{k}\left(  x\right)  \right)  ^{2}\dfrac
{1-x}{1+x}dx+\beta\left(  \sum_{k=0}^{n}a_{k}P_{k}\left(  1\right)  -1\right)
\]
where $\beta$ is the Lagrangian multiplier, \cite{Reh}, \cite{Szeg} we have%
\begin{equation}
\frac{\partial\varphi}{\partial a_{k}}=2\int\limits_{-1}^{1}a_{k}P_{k}%
^{2}\left(  x\right)  \dfrac{1-x}{1+x}dx+\beta P_{k}\left(  1\right)  =0
\label{Lagrange}%
\end{equation}
and%
\begin{equation}
\sum_{k=0}^{n}a_{k}P_{k}\left(  1\right)  =1 \label{conditionn}%
\end{equation}
Denoting by%
\[
\left\Vert P_{n}\right\Vert ^{2}=%
{\displaystyle\int\limits_{-1}^{1}}
P_{n}^{2}\left(  x\right)  \dfrac{1-x}{1+x}dx=\frac{2\left(  n+1\right)  ^{2}%
}{n\left(  n-1\right)  \left(  2n-1\right)  }\text{ \ \ \ \ \ ,}%
n=2,3,4........
\]
we easily find,by \eqref{conditionn}
\[
a_{k}=\frac{P_{k}\left(  1\right)  }{\left\Vert P_{k}\right\Vert ^{2}}\frac
{1}{%
{\displaystyle\sum\limits_{j=2}^{n}}
\frac{P_{j}\left(  1\right)  ^{2}}{\left\Vert P_{j}\right\Vert ^{2}}}%
\]
so that the minimum value $M$ of the integral \eqref{ExtrPbm} under the
aforementioned constraint is%
\begin{equation}
M=\int\limits_{-1}^{1}\left(  \sum_{k=2}^{n}\frac{P_{k}\left(  1\right)
M}{\left\Vert P_{k}\right\Vert ^{2}}P_{k}\left(  x\right)  \right)  ^{2}%
\dfrac{1-x}{1+x}dx=\frac{1}{%
{\displaystyle\sum\limits_{j=2}^{n}}
\frac{P_{j}\left(  1\right)  ^{2}}{\left\Vert P_{j}\right\Vert ^{2}}}
\label{Mvalue}%
\end{equation}
and%
\begin{equation}
f_{n}\left(  x\right)  =\frac{1}{%
{\displaystyle\sum\limits_{j=2}^{n}}
\frac{P_{j}\left(  1\right)  ^{2}}{\left\Vert P_{j}\right\Vert ^{2}}}%
\sum_{k=2}^{n}\frac{P_{k}\left(  1\right)  P_{k}\left(  x\right)  }{\left\Vert
P_{k}\right\Vert ^{2}} \label{Solution}%
\end{equation}
becomes to the following solution of above extremal problem :
\begin{equation}
f_{n}\left(  x\right)  =\sum_{k=2}^{n}\frac{MP_{k}\left(  1\right)
P_{k}\left(  x\right)  }{\left\Vert P_{k}\right\Vert ^{2}} \label{Solution2}%
\end{equation}

\begin{theorem}
the integral
\begin{equation}
I_{n}=\int\limits_{-1}^{1}\left(  F_{n}\left(  x\right)  \right)  ^{2}%
\dfrac{1-x}{1+x}dx \label{Integral}%
\end{equation}

where $F_{n}\left(  x\right)  $ is any real polynomial of degree $n$ such that
$F_{n}(1)=1,$reaches its minimum value
\begin{equation}
M=\frac{2}{%
{\displaystyle\sum\limits_{j=2}^{n}}
j\left(  j-1\right)  \left(  2j-1\right)  } \label{Valuem}%
\end{equation}

if and only if%
\begin{equation}
F_{n}\left(  x\right)  =\frac{2}{%
{\displaystyle\sum\limits_{j=2}^{n}}
j\left(  j-1\right)  \left(  2j-1\right)  }\sum_{k=2}^{n}\frac{k\left(
k-1\right)  \left(  2k-1\right)  }{2\left(  k+1\right)  }P_{k}\left(
x\right)  \label{Fnnx}%
\end{equation}

$\left\{  F_{n}\left(  x\right)  \right\}  _{n=2,3,4......}$are orthogonal
over $\left[  -1\text{ \ \ \ }1\right]  ,$with respect to the weight function
$x\longrightarrow-\dfrac{\left(  x-1\right)  ^{2}}{1+x}.$Hence%
\begin{equation}
M=\frac{1}{K_{n}\left(  0,0\right)  } \label{Kernelm}%
\end{equation}

and%
\begin{equation}
F_{n}\left(  x\right)  =\frac{K_{n}\left(  x,0\right)  }{K_{n}\left(
0,0\right)  } \label{Kernelf}%
\end{equation}

\end{theorem}

\begin{proof}
Using \eqref{ExtrPbm}, \eqref{Mvalue}, \eqref{Solution},
\eqref{Pnat1},\eqref{NormQn},gives the minimum value%
\[
M=\frac{2}{%
{\displaystyle\sum\limits_{j=2}^{n}}
j\left(  j-1\right)  \left(  2j-1\right)  }%
\]

and%
\[
F_{n}\left(  x\right)  =\frac{2}{%
{\displaystyle\sum\limits_{j=2}^{n}}
j\left(  j-1\right)  \left(  2j-1\right)  }\sum_{k=2}^{n}\frac{k\left(
k-1\right)  \left(  2k-1\right)  }{2\left(  k+1\right)  }P_{k}\left(
x\right)
\]

and this completes the proof of Theorem .
\end{proof}

\begin{theorem}
Let $f$ be an increasing function on $\left[  -1\text{ \ \ \ }1\right]  $
,with $f\left(  a\right)  =-1$ and $f\left(  b\right)  =1,$such that $a<b$ and
$\varphi$ a nonnegative weight function on the same interval, such that the
integral
\[
\int\limits_{-1}^{1}f\left(  x\right)  ^{n}\varphi\left(  x\right)  dx\text{
\ \ \ \ \ \ \ \ \ \ \ }\left(  n\geq0\right)
\]

exists; Then the sequence of functions $x\mapsto P_{0}\left(  f\left(
x\right)  \right)  ,x\mapsto P_{1}\left(  f\left(  x\right)  \right)
,...x\mapsto P_{n}\left(  f\left(  x\right)  \right)  ...$that minimizes the
integrals%
\begin{equation}
I_{n}=\int\limits_{a}^{b}q_{n}\left(  f\left(  x\right)  \right)  ^{2}%
\varphi\left(  x\right)  dx \label{In}%
\end{equation}

for all polynomial : $q_{n}\left(  x\right)  =b_{0}+b_{1}x+......b_{n}x^{n}$,
forms an orthogonal system on $\left[  a\text{ \ \ \ }b\right]  $ in respect
of \ $\varphi.$Where%
\begin{equation}
\frac{\varphi\left(  x\right)  }{f^{\prime}\left(  x\right)  }=\dfrac
{1+f\left(  x\right)  }{1-f\left(  x\right)  } \label{weight}%
\end{equation}

i.e,%
\[
\int\limits_{a}^{b}P_{n}\left(  f\left(  x\right)  \right)  P_{m}\left(
f\left(  x\right)  \right)  \dfrac{1+f\left(  x\right)  }{1-f\left(  x\right)
}f^{\prime}\left(  x\right)  dx=0\text{ \ \ \ \ \ \ \ \ ,}n=0,1,2.....\left(
n\neq m\right)
\]

If
\begin{equation}
f\left(  x\right)  =\frac{4x^{3}}{\left(  x^{2}+1\right)  ^{2}} \label{fff}%
\end{equation}

satisfie $f\left(  -1\right)  =-1,f\left(  1\right)  =1,$then%
\[
\int\limits_{-1}^{+1}P_{n}\left(  \frac{4x^{3}}{\left(  x^{2}+1\right)  ^{2}%
}\right)  P_{m}\left(  \frac{4x^{3}}{\left(  x^{2}+1\right)  ^{2}}\right)
\varphi\left(  x\right)  dx=0\text{ \ \ \ \ \ \ \ \ ,}n=0,1,2.....\left(
n\neq m\right)
\]

where%

\begin{equation}
\varphi\left(  x\right)  =\frac{\left(  x^{2}+1\right)  ^{2}+4x^{3}}{\left(
x^{2}+1\right)  ^{2}-4x^{3}}\left(  \frac{12x^{2}}{\left(  x^{2}+1\right)
^{2}}-\frac{4x}{\left(  x^{2}+1\right)  ^{3}}\right)  \label{wffff}%
\end{equation}

\end{theorem}

\begin{proof}
the polar Legendre polynomials $\left\{  P_{n}\right\}  _{n=0,1,2,...}$ are
orthogonal on $[-1$ $\ \ \ 1]$ in respect of
\[
t\mapsto\psi\left(  t\right)  =\dfrac{1+t}{1-t}%
\]

i-e%
\[
\int\limits_{-1}^{1}P_{n}\left(  t\right)  P_{m}\left(  t\right)  \dfrac
{1+t}{1-t}dt=0\text{ \ \ \ \ \ \ \ \ ,}n,m=0,1,2.....\left(  n\neq m\right)
\]

Substituting $f(x)=t$ \ in \eqref{In} we have%
\begin{equation}
I_{n}=\int\limits_{-1}^{1}q_{n}\left(  t\right)  ^{2}\frac{\varphi\left(
f^{-1}\left(  t\right)  \right)  }{f^{\prime}\left(  f^{-1}\left(  t\right)
\right)  }dt\text{ \ \ \ \ \ \ ,}n=0,1,2....., \label{eee}%
\end{equation}

if%
\[
\frac{\varphi\left(  f^{-1}\left(  t\right)  \right)  }{f^{\prime}\left(
f^{-1}\left(  t\right)  \right)  }=\dfrac{1+t}{1-t}%
\]

Now comming back to the old variable with according to Theorem 1, the
minimizing functions
\[
x\mapsto P_{0}\left(  f\left(  x\right)  \right)  ,x\mapsto P_{1}\left(
f\left(  x\right)  \right)  ,x\mapsto P_{2}\left(  f\left(  x\right)  \right)
,...x\mapsto P_{n}\left(  f\left(  x\right)  \right)  ....
\]

that minimize \eqref{In} form an orthogonal system on $\left[  a\text{
\ \ \ }b\right]  $ in respect of \ $\varphi.$Therefore we denote as%
\[
x\longrightarrow P_{k}\left(  \frac{4x^{3}}{\left(  x^{2}+1\right)  ^{2}%
}\right)  ,\text{ \ \ }k=2,3,4........
\]

form an orthogonal system on $\left[  -1\text{ \ \ \ }1\right]  $ in respect
of \ $\varphi.$%
\[
\varphi\left(  x\right)  =\frac{\left(  x^{2}+1\right)  ^{2}+4x^{3}}{\left(
x^{2}+1\right)  ^{2}-4x^{3}}\left(  \frac{12x^{2}}{\left(  x^{2}+1\right)
^{2}}-\frac{4x}{\left(  x^{2}+1\right)  ^{3}}\right)
\]
i-e%
\[
\int\limits_{-1}^{+1}P_{n}\left(  \frac{4x^{3}}{\left(  x^{2}+1\right)  ^{2}%
}\right)  P_{m}\left(  \frac{4x^{3}}{\left(  x^{2}+1\right)  ^{2}}\right)
\varphi\left(  x\right)  dx=0\text{ \ \ \ \ \ \ \ \ ,}n=0,1,2.....\left(
n\neq m\right)
\]

and this completes the proof of Theorem
\end{proof}

\begin{example}
Let $f$ be an increasing function on $\left[  -1\text{ \ \ \ }1\right]  $
,with $f\left(  u\right)  =-1$ and $f\left(  v\right)  =1$ and $\varphi$ a
nonnegative weight function on the same interval, such that
\begin{equation}
f\left(  x\right)  =\frac{ax+b}{cx+d}\text{ \ \ \ \ \ , }x\neq-\frac{d}{c}
\label{fgfg}%
\end{equation}

then%
\[
x\mapsto P_{0}\left(  f\left(  x\right)  \right)  ,x\mapsto P_{1}\left(
f\left(  x\right)  \right)  ,x\mapsto P_{2}\left(  f\left(  x\right)  \right)
,...x\mapsto P_{n}\left(  f\left(  x\right)  \right)  ....
\]

form an orthogonal system on $\left[  u\text{ \ \ \ }v\right]  $ in respect of
\ $\varphi$.where%
\[
\varphi\left(  x\right)  =\frac{ad-bc}{\left(  cx+d\right)  ^{2}}%
\dfrac{\left(  a+c\right)  x+b+d}{\left(  c-a\right)  x-b+d}%
\]

i-e%
\[
\int\limits_{u}^{v}P_{n}\left(  \frac{ax+b}{cx+d}\right)  P_{m}\left(
\frac{ax+b}{cx+d}\right)  \frac{ad-bc}{\left(  cx+d\right)  ^{2}}%
\dfrac{\left(  a+c\right)  x+b+d}{\left(  c-a\right)  x-b+d}dx=0\text{
\ \ \ \ \ },\left(  n\neq m\right)
\]

\end{example}

\end{document}